\parskip=4pt
\mathsurround=2pt
\input amssym.def
\input amssym

\def\v{{\bf v}}
\def\w{{\bf w}}
\def\O{{\bf o}}
\def\R{{\bf R}}
\def\sgn{{\rm sgn}}
\def\QED{~\rlap{$\sqcup$}$\sqcap$ \smallskip}
\def\QP{\smallskip\leftskip=.4in\rightskip=.4in\noindent}
\def\ref{\hangindent=1pc \hangafter=1 \noindent}
\def\s{\thinspace}

\font\bit=cmssi12 at 12 truept
\font\sma=cmr10 at 10 truept

\input psfig

\centerline{\bf A Monotonicity Conjecture for Real Cubic Maps.\footnote{$^1$}
{\sma Based on lectures by Milnor
at the NATO Advanced Study Institute on Real and Complex
Dynamical Systems, Hiller{\o}d, June 1993.}
}\medskip

\centerline{Silvina P. Dawson, Roza Galeeva, John Milnor\footnote{$^2$}
{\sma Partially supported by the Miller Institute of the University
of California at Berkeley during the preparation of this paper.}
and Charles Tresser}

\bigskip

\centerline{\bf 1. Introduction.}\smallskip

This will be an outline of work in progress. We study the
conjecture that the topological entropy of a real cubic map depends
``monotonely'' on its parameters, in the sense that each locus of
constant entropy in parameter space is a connected set.

Section 2 sets the stage by describing the parameter triangle $\overline T$ for real
cubic maps, either of shape $+-+$ or of shape $-+-$, and by describing basic
properties of topological entropy. Section 3 describes the
monotonicity problem for the topological entropy function, and states the
Monotonicity Conjecture. Section 4 describes the family of stunted
sawtooth maps, and proves the analogous conjecture for this family. Sections
5 begins to relate these two families by describing the `bone' structure
in the parameter triangle. By definition, a {\bit bone\/} $B_\pm(\O)$ in the
triangle $\overline T$ is the set of parameter points $\v$ such that a specified critical
point (left or right) of the associated bimodal map belongs to a
periodic orbit with specified order type $\O$. (Compare  [Ma\s T].) It is
conjectured
that every bone is a simple connected arc in $\overline T$. Although we cannot
prove either of these conjectures for cubic maps, we do show that\smallskip

\centerline{\bit Generic Hyperbolicity~ $\Rightarrow$~ Connected Bone
Conjecture~ $\Rightarrow$~ Monotonicity Conjecture}

\noindent (see Theorems 3 and 4 in Sections 7, 8). The paper concludes with
a brief Appendix on computation.

This material will be presented in more detail in a later paper [DGMT].
\bigskip

\centerline{\bf 2. The Parameter Triangle $\overline T$ and the Topological
Entropy Function.}\smallskip

Let $f$ be a cubic map of the unit interval $I=[0,1]$. We will
always assume that $f$ maps the boundary of $I$ into itself.
To fix our ideas, we consider only those maps which have
{\bit shape\/} $+-+$; that is, $f$ must first increase, then decrease, and
then increase.
Thus the leading coefficient must be positive, $f$ must have
critical points $c_1<c_2$ in the interior of $I$, and
both boundary points must be fixed by $f$. (All of the discussion
which follows could easily be modified so as to apply also
to maps of shape  $-+-$; these are dynamically quite different, since
the two boundary points must form a period two orbit.) For maps of shape
$+-+$, evidently
the corresponding {\bit critical values} $v_i=f(c_i)$ must satisfy
$$1\ge v_1>v_2\ge 0~. \eqno (1)$$
\smallskip 

{\QP{\bf Lemma 1.} {\it Given any pair  ${\bf v}=(v_1,v_2)$ satisfying the
inequalities $(1)$, there is one and only one cubic map which fixes the
boundary of $I$ and has critical values $(v_1,v_2)$.}\smallskip}

{\bf Proof Outline.} It is easy to check that any
real cubic map with distinct real critical points can be written uniquely as
$$	f(x)~=~ c\,F(ax+b)+d $$
with $a>0$, where
$$	F(\xi) ~=~3\,\xi^2-2\,\xi^3$$
is the unique cubic map with
fixed critical points $0$ and $1$. Note that $f$ has critical values
$v_1=d$ and $v_2=c+d$. We can solve these linear equations for $c$ and $d$,
and then solve the required cubic equations
$f(0)=0$ and $f(1)=1$ for $b$ and $a+b$.\QED

{\bf Definitions.} This cubic map with {\bit critical value vector}
${\bf v}=(v_1,v_2)$ will be denoted by $f=f_{\bf v}$. It is often convenient
to allow the limiting case $v_1=v_2$ also. This corresponds to allowing
degenerate cubic maps, for which $c_1=c_2$. The
compact set $\overline T\subset\R^2$ consisting of all pairs $(v_1\,,\,v_2)$ with
$1\ge v_1\ge v_2\ge 0$ will be called
the {\bit parameter triangle\/} for cubic maps of shape $+-+$.

\midinsert
\centerline{\psfig{figure=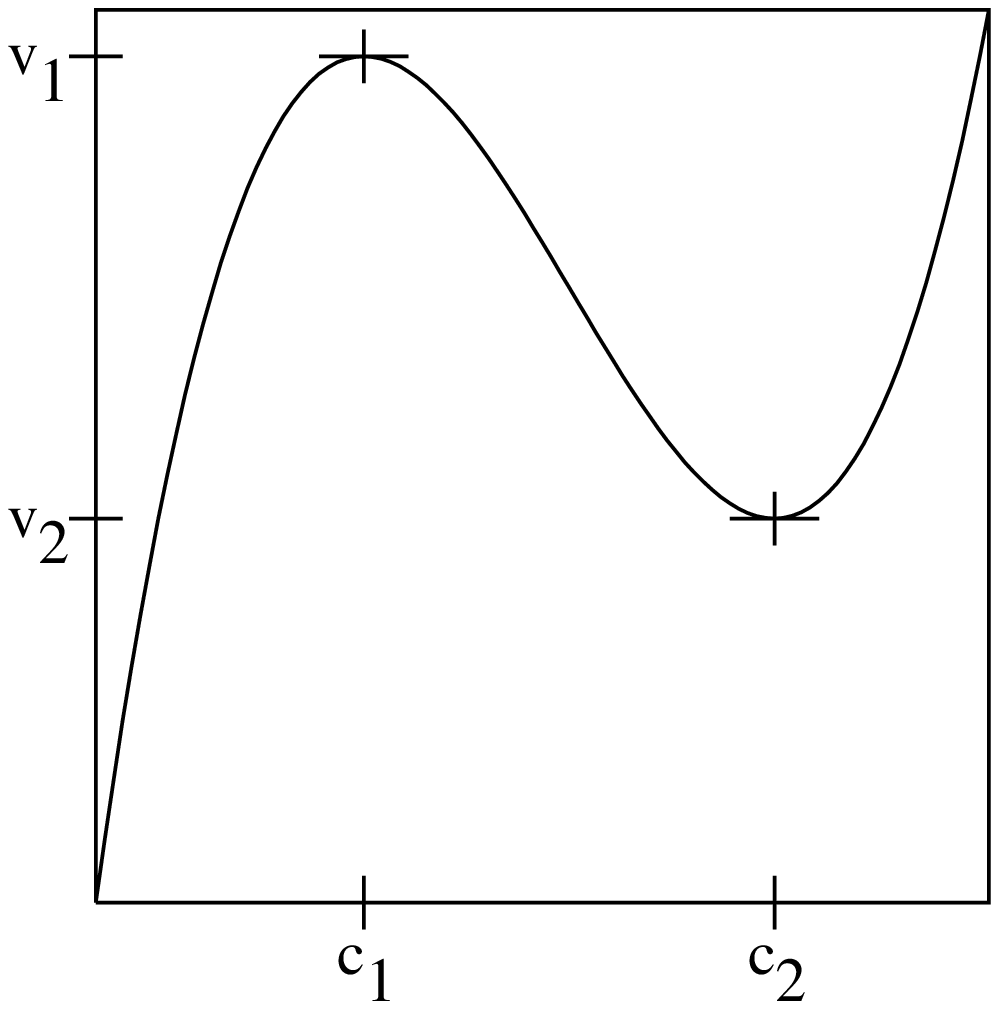,height=2in}
\psfig{figure=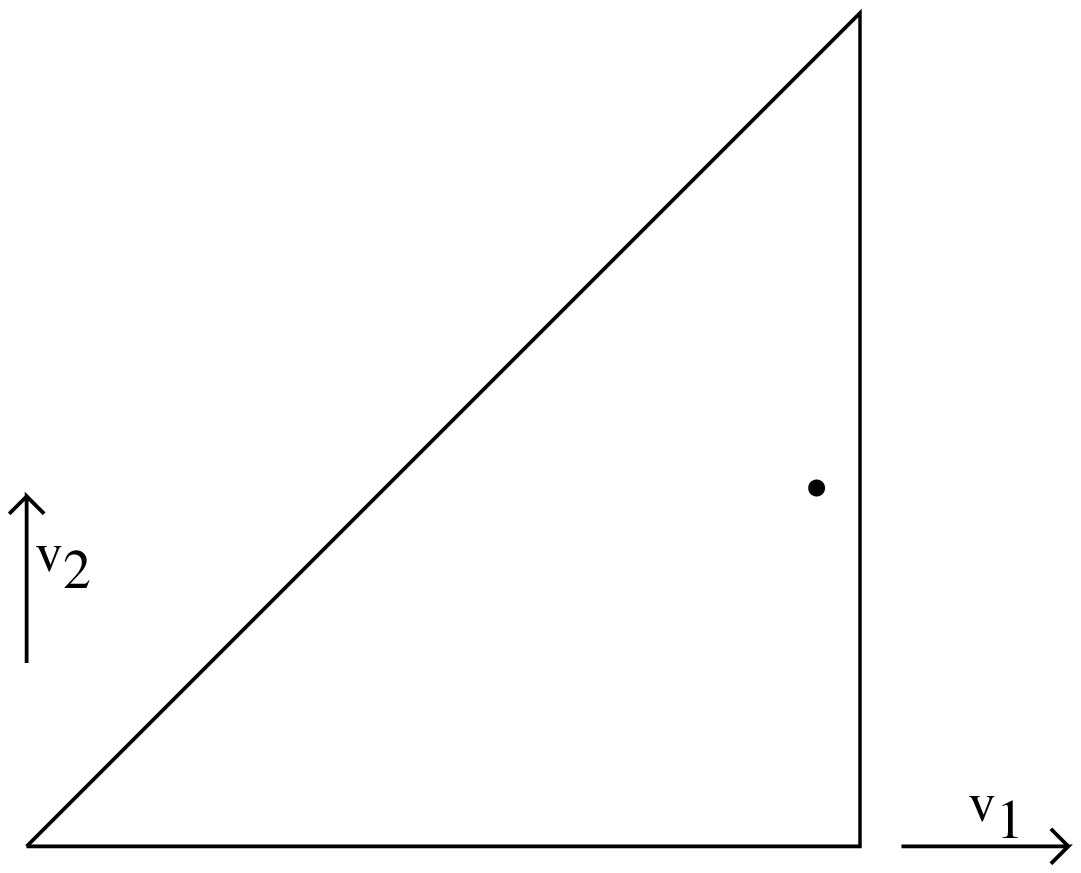,height=2in}}
\centerline{\bit A cubic map $f_{\bf v}$, and the corresponding point
$\v=(v_1,v_2)$ in parameter space.}
\endinsert

{\bf Remarks.} The analogue of Lemma 1 for cubic maps of shape $-+-$
can be proved by essentially the same argument. Corresponding
statements for higher degree polynomials with distinct real critical points
are also true.
For example, such a polynomial can be constructed uniquely, from its
critical value vector,
by constructing the Riemann surface of the corresponding complex
polynomial map. A purely
real proof may be found in [dM\s vS, p. 120]. (A somewhat
simpler real proof has been given by Douady and Sentenac, unpublished.)
For further information, see [DGMT].
\medskip

How can we measure the dynamic ``complexity'' of a map $f_\v:I\to I$,
and how does this complexity vary as the critical value
vector $\bf v$ varies within the parameter
triangle $\overline T$? One measure of complexity
would be the numbers of periodic points of various
periods. A particularly useful measure of complexity is provided by
the topological entropy $h$. For our purposes, the {\bit topological
entropy\/} of a piecewise-monotone map can be defined by the formula
$$ h(f)~=~\lim_{n\to\infty}{\log\ell(f^{\circ n})\over n}~,\eqno (2) $$
where $\ell(f^{\circ n})$ is the number of {\bit laps\/} of the
$n\!$-fold iterate, that is the number of maximal intervals of monotonicity.
(Compare [Ro], [M\s Sz]. For computation of $h$, see [B\s K], [BST].)
For non-linear polynomial maps, or more generally for piecewise-monotone
maps with at most finitely many non-repelling periodic orbits, $h(f)$
can be identified with the number
$$	h_{\rm per}(f)~=~\limsup_{n\to\infty}{\log\#{\rm fix}(f^{\circ n})
\over n}~,\eqno (3) $$
where $\#{\rm fix}$ is the number of fixed points.\footnote{$^1$}
{\sma See [M\s Sz], [M\s Th], as well as [dMvS, p.268].
It is possible that the equation $h=h_{\rm per}$ is true for a
$C^r\!$-generic map in any dimension, but no proof is known. (Compare
 [B, p. 23].)}
(Compare Lemma 7 in \S8.)
The entropy varies continuously under bimodal $C^1\!$-deformation.\footnote
{$^2$}{\sma See [M\s Sz], [M\s Th], [dM\s vS]. Conjecturally, entropy
remains continuous under $C^1\!$-deformation as long as the number of critical
points remains bounded; but even in the $C^r\!$-case it
definitely can jump discontinuously if the number of critical points is
unbounded and if $r<\infty$. See [M\s Sz].
For maps in dimension $\ge 2$ or for diffeomorphisms in dimension $\ge 3$,
the entropy can also drop discontinuously, even in the
$C^\infty$ case. See [K], [Mis], [N], [Y] for further information.}

In practice, it is often more convenient to work with the quantity
$$	s~=~\exp(h)~=~\lim_{n\to\infty} \root n \of {\ell(f^{\circ n})}~, $$
sometimes known as the {\bit growth number\/} of $f$.
For an $m\!$-{\bit modal\/} map, that is for a map with $m+1$ laps,
this number $s$ lies in the closed interval
$[1\,,\, m+1]$. In the special case
of a piecewise linear map with $~|{\rm slope}|={\rm constant}\ge 1$, the
growth number $s$ is precisely equal to this constant $|{\rm slope}|$.
\bigskip

\centerline{\bf 3. The Monotonicity Problem.}\smallskip

In the quadratic case, it is known that the number of period $p$ points
for an interval map $x\mapsto 4vx(1-x)$ increases monotonically as the
critical value
parameter $v\in [0,1]$ increases. (Proofs of this result have been given by
Sullivan, Douady and Hubbard, and by Milnor and Thurston.
Compare [DH2, n$^{\!\rm o}$VI],
[M\s Th], [D], as well as [dM\s vS].) Hence the entropy $h$
also increases monotonically with $v$. Compare the picture below. 


\midinsert
\centerline{\qquad\qquad\psfig{figure=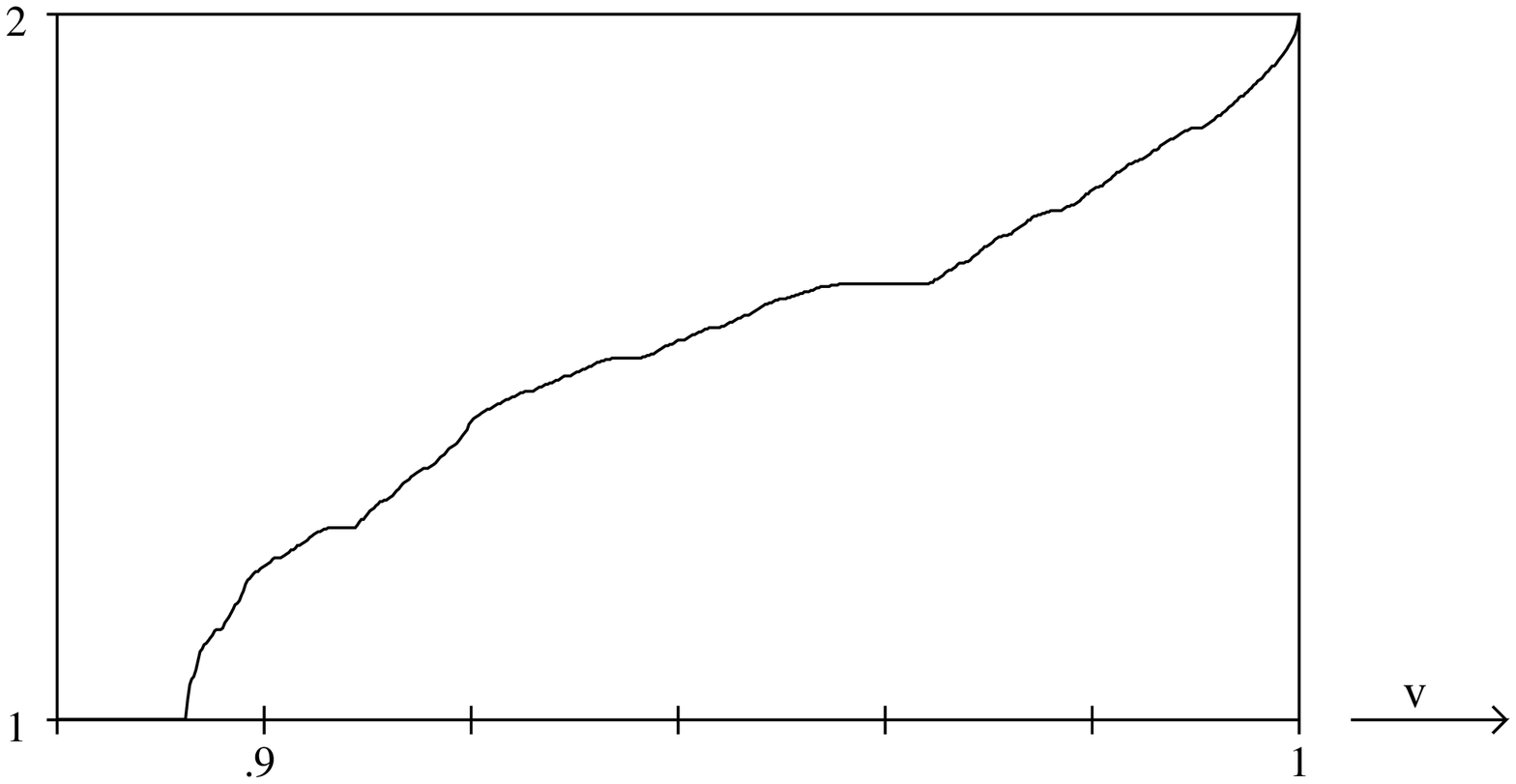,height=2.3in}}
\centerline{\bit Graph of $\,s=\exp(h)\,$ as a function of $\,v\,$
for the family of maps $\,x\mapsto 4vx(1-x)$.}
\endinsert
\eject

Is there some analogous statement for the two-parameter family of cubic
maps $f_\v$? Can we find curves through an arbitrary point of $\overline T$
along which
the complexity increases monotonically? (Compare [D\s G], [DGK], [DGKKY],
[DGMT].) Is there some sense
in which the topological entropy function
$$	{\bf v}~\mapsto~h(f_{\bf v}) $$
is a ``monotone'' function on the parameter triangle $\overline T$? To formulate this
question more precisely, we make the following definition. Fix some
constant $h_0$ in the closed interval $[0,\log 3]$. By the {\bit
$h_0\!$-isentrope\/} for the family of cubic maps $f_\v$ we will mean the
set consisting of all parameter values $\v\in \overline T$ for which the
topological entropy $h(f_\v)$ is equal to $h_0$. (In the illustrations,
it will be convenient to work with $s=\exp(h)$ rather than $h$.)
Since the entropy function
is continuous, note that each isentrope is a compact subset of $\overline T$.

{\QP{\bf Monotonicity Conjecture} for the family of real cubic maps
$f_\v$. {\it Every isentrope
$	\{\v\in T~:~h(f_\v)=h_0\,\} $
for this family is a connected set.}\bigskip}

\midinsert
\centerline{\psfig{figure=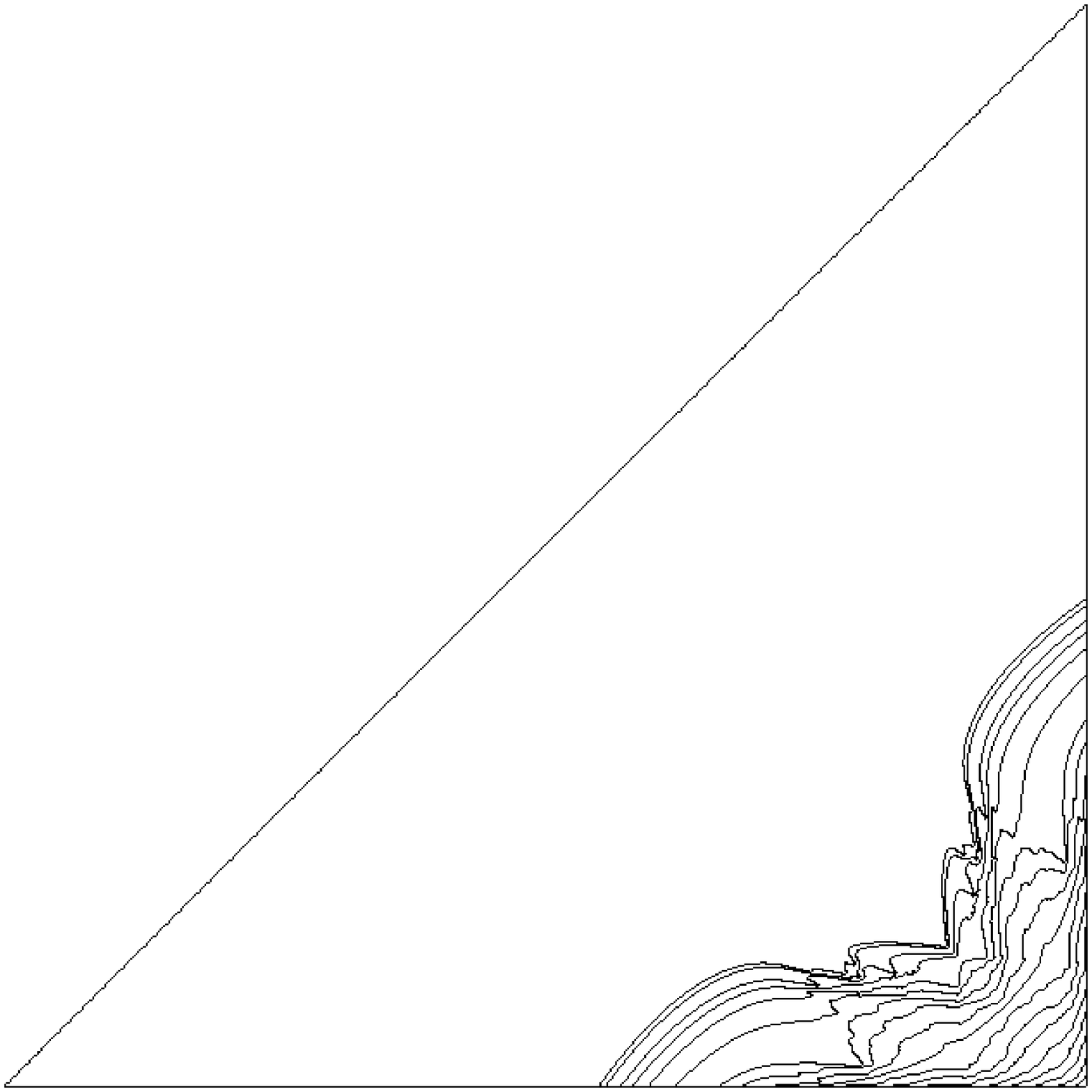,height=3.2in}}
\centerline{\bit Isentropes $\,s={\rm constant}\,$ in the parameter
triangle for cubic maps of shape $+-+$.}

\centerline{\bit (Contour interval: $\Delta s=0.1$. Visible isentropes:
$s=1.1,\,1.2,\,1.3,\,\ldots\,,\,2.8$.)}
\endinsert

\midinsert
\centerline{\psfig{figure=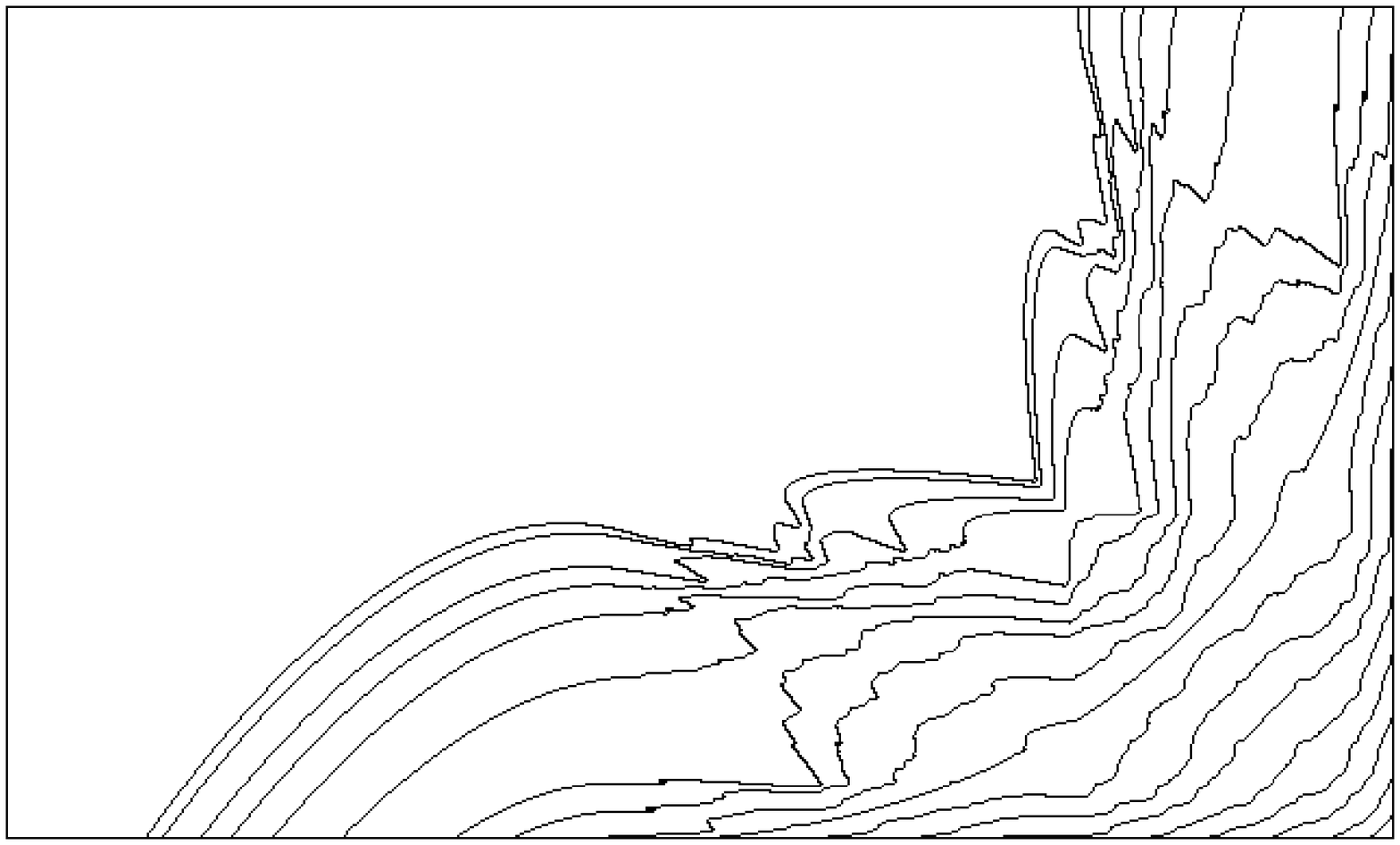,height=2.2in}}

\vskip 5pt plus 2 ex
\centerline{\psfig{figure=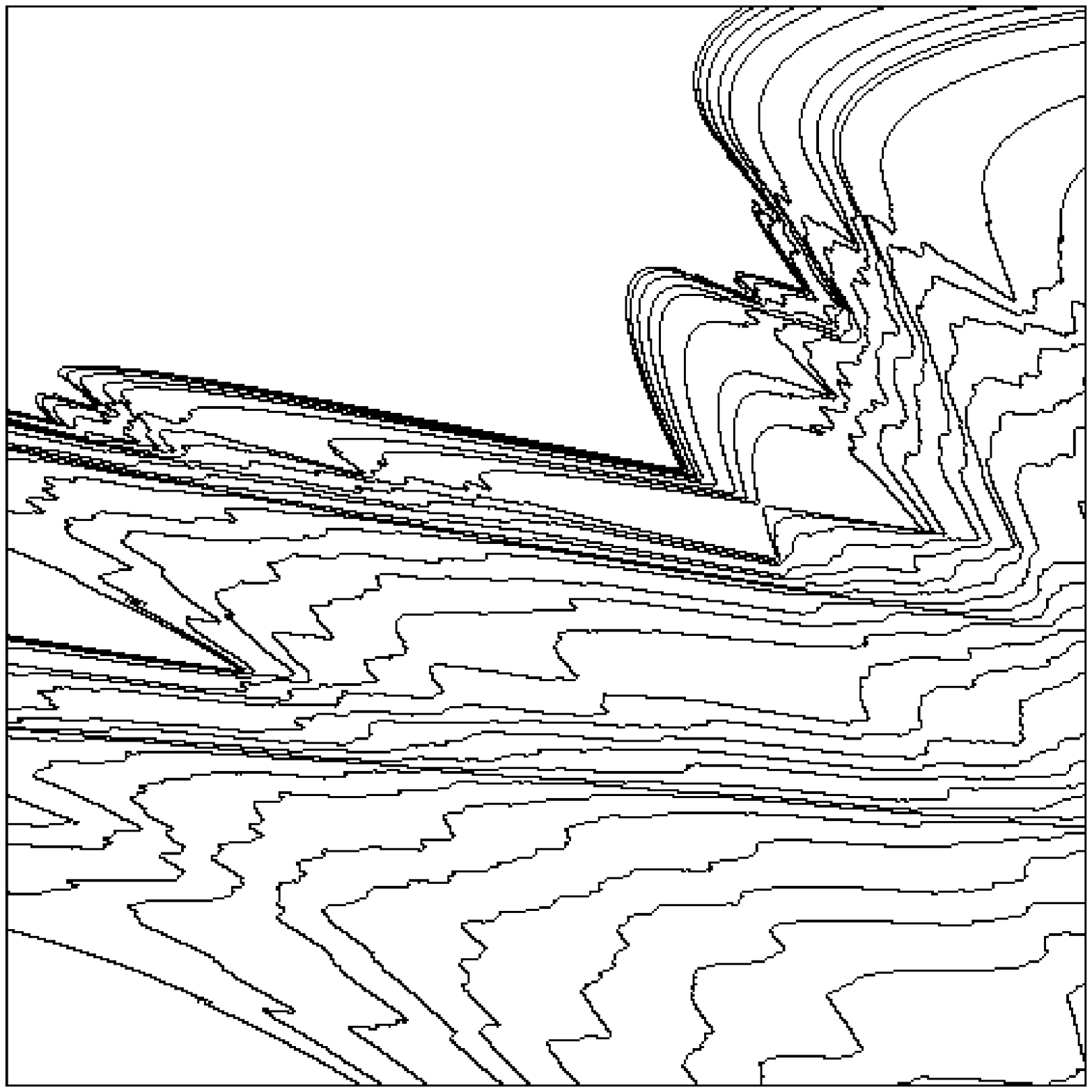,height=3.5in}}

{\QP\bit Magnified portions of the previous figure. Above:
The lower right region\break $[.5,1]\times[0,.3]$, again with contour interval
$\Delta s= 0.1$. Below: Detail showing the region $[.74,.8]\times[.07,.13]$
near the center of this picture, with
$\Delta s=.02$.\par}
\endinsert

(Compare \S6, as well as [M1, p.13].)
Evidently this conjecture describes a weak form of monotonicity for this
two-parameter family. Such a
connected isentrope could be a simple arc with endpoints on
the boundary of T, or perhaps could have a more complicated non-locally
connected topology although
this has not been observed. It is also certainly possible for it to be a
compact set with interior points.
In the limiting case $h_0=0$, the isentrope is the large white region
in the picture above, containing
the entire upper left hand edge $v_1=v_2$. (Compare [Ma\s T].)
In the other limiting case
$h_0=\log 3$, it reduces to the single corner point $\v=(1,0)$. If the
Conjecture is true, then for each
$0<h_0<\log 3$ the $h_0\!$-isentrope must cut $\overline T$ into
two connected pieces, one with $h<h_0$ and one with $h>h_0$.
\smallskip

Another interesting consequence would be a ``maximum and minimum
principle'' for entropy:
{\it If the conjecture is true, then
the maximum and minimum values for the entropy function on any closed
region $U\subset \overline T$ must occur on
the boundary $\partial U$.\/} In fact every
value of entropy which occurs in $U$ must occur already on $\partial U$.
This follows since, by continuity, every
value of entropy between $0$ and $\log 3$
must occur on $\partial \overline T$. Evidently a connected
isentrope which contains points both inside and outside $U$ must also
intersect $\partial U$.

These questions are quite difficult. As a step towards understanding, we
will first consider a different
family of maps for which they are very much easier.
\bigskip\bigskip 

\centerline{\bf 4. The Stunted Sawtooth Family.}\smallskip

By the {\bit sawtooth map} of shape $+-+$ we mean the unique map $S:I\to I$
which is piecewise linear with slope alternately $+3\,,\,-3$ and $+3$.
This is a bimodal map for which the topological entropy takes the
largest possible value $h=\log 3$. {\bf Definition:}
Given any critical value vector ${\bf w}=(w_1,w_2)$ satisfying the usual
inequalities $1\ge w_1>w_2\ge 0$, we obtain the {\bit stunted sawtooth map}
$S_{\bf w}$ from $S$
by cutting off the top and bottom at heights $w_1$ and $w_2$, as illustrated
below. (Compare [Gu].)
As with the cubic family, it is often convenient to allow the limiting
case $w_1=w_2$. For the stunted sawtooth family,
this means that we may allow the
two horizontal plateaus to come together.

\midinsert
\centerline{\psfig{figure=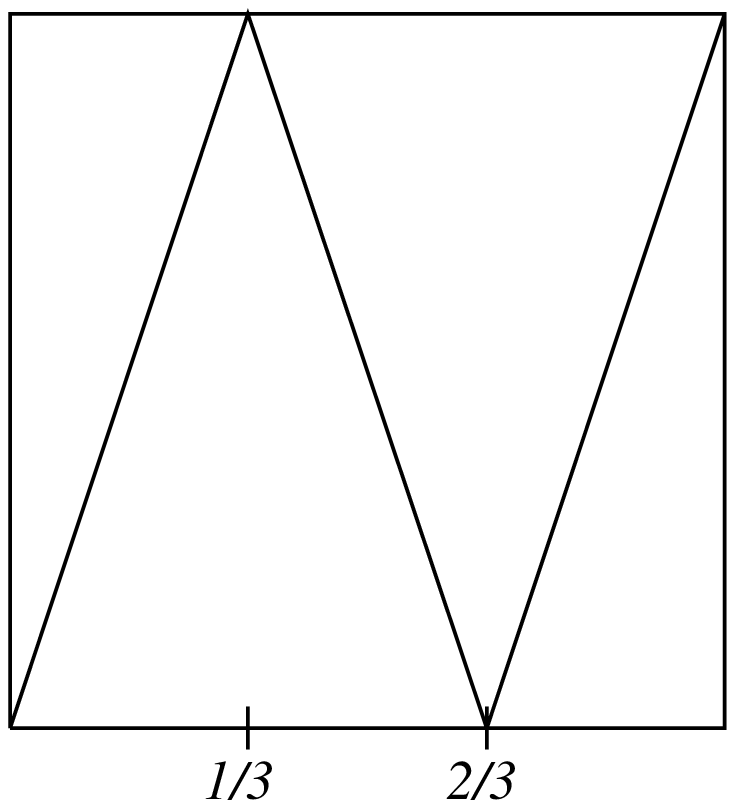,height=2.2in}
\psfig{figure=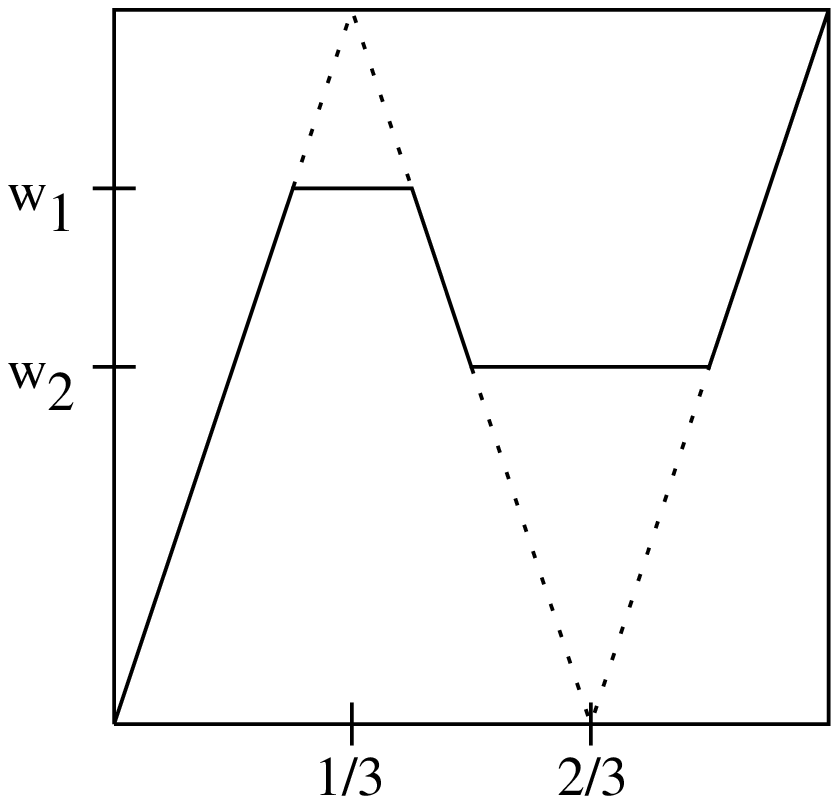,height=2.2in}}
\centerline{\bit The sawtooth map  $S$\qquad\qquad A stunted sawtooth map
$S_{\bf w}$.}
\endinsert

For this family it is easy to see that any increase in the parameter 
$w_1$, or any increase in $1-w_2$, can only increase the
complexity of the mapping. (Compare [BMT], [Ga].)
For example, as we increase $w_1$ with
fixed $w_2$, no periodic orbit can disappear:
If a given periodic orbit misses the left hand
plateau, then it remains unchanged as we increase $w_1$, while if it
hits this plateau then it deforms continuously as we increase $w_1$.
Similarly, the topological entropy can never decrease as we increase
$w_1$ or $1-w_2$. It will be convenient to define a simple
partial ordering for the parameter triangle $\overline T$ as follows:
$$	 \w\ll\w'\quad\Longleftrightarrow\quad
  w_1\le w'_1\quad{\rm and}\quad 1-w_2\le 1-w'_2~.\leqno {\bf\qquad
  Definition:}$$
Then it follows from the discussion above that
$$	\w\ll\w'\quad\Longrightarrow\quad h(S_\w)\le h(S_{\w'})~. \eqno (4) $$
\eject

{\bf Remark.} Let us temporarily extend the discussion to more general
piecewise-mono\-tone maps. By the {\bit shape\/} of an $m\!$-modal map
we mean an alternating
sequence of $m+1$ signs, starting with either $+$ or $-$ according as
the map is increasing or decreasing on its initial lap.
The above construction for bimodal maps of shape $+-+$ extends easily to
$m\!$-modal maps for any $m\ge 1$ and for either one of the two possible
shapes.

The stunted sawtooth
family is very closely related to kneading theory. To make this precise,
we will need the following. Again consider an $m\!$-modal map
with any $m\ge 1$ and with either one of the two possible $m\!$-modal
shapes.

{\bf Definition.} By the {\bit kneading data\/} associated with an $m\!$-modal
map $f$ we will mean its shape, together with the collection of signs
$$    \sgn\big(f^{\circ n}(c_i) -c_j\big)~~\in~~\{-1,0,1\} $$
for $n>0$ and $1\le i,j\le m$, where the $c_i$ are the critical points of $f$.

To extend this definition to the case of an $m\!$-modal
stunted sawtooth map, we simply define the ``critical points'' to be the
center points $\hat c_i=i/(m+1)~$ of the plateaus, for\break
$1\le i\le m$. With this
definition, we can make the following assertion.

{\QP{\bf Lemma 2.} \it To any $m\!$-modal map $f$ there is associated a canonical
stunted
sawtooth map $S_\w$ which has exactly the same kneading data.\smallskip}

The proof can be outlined as follows (details in [DGMT]). Let $S$ be the
$m\!$-modal sawtooth map with the same shape, with critical points $\hat c_i=
i/(m+1)$. First consider a point $x$ in
the domain of definition of $f$ which is not pre-critical. That is, we
assume that the orbit $\{x\,,\,f(x)\,,\,f^{\circ 2}(x)\,,\,\ldots\}$ does
not contain any critical point. Then
there exists one and only one point $\hat x\in[0,1]$ so that
the {\bit itinerary\/} of $\hat x$ under $S$ is the same as the itinerary
of $x$ under $f$. By definition, this means that
$$	\sgn\big(f^{\circ n}(x)-c_j\big) ~~=~~
   \sgn\big(S^{\circ n}(\hat x)-\hat c_j\big)\eqno (5) $$
for every $n\ge 0$ and for every critical point $c_j$ of $f$. In the case of
a pre-critical point $x$, we must weaken this condition slightly
by requiring equation (5) only up to the first $n$ for which $f^{\circ n}(x)$
is critical. Then again there is a unique associated $\hat x$.
Now let $v_1\,,\,\ldots\,,\, v_m$ be the critical values of $f$. Then the
associated points $w_i=\hat v_i$ are the critical values for the required
stunted sawtooth map $S_\w$.\QED

In particular, for each $n>0$ we have a matrix equality
$$	\big[ \sgn\big(f^{\circ n}(c_i)-c_j\big)\big]~=~\big[\sgn\big(
  S_\w^{\circ n}(\hat c_i)-\hat c_j\big)\big]~, $$
where the $c_i$ are the critical points of $f$ and where the $\hat
c_i=i/(m+1)$ are the critical points of $S_\w$.\medskip
\eject

Now let us again specialize to $+-+$ bimodal maps. We will show that
the Monotonicity Conjecture for the stunted sawtooth family is true:

{\QP{\bf Theorem 1.} {\it For every constant $0\le h_0\le\log 3$, the
isentrope\vskip -.1in
$$	{\cal I}(h_0)~=~\{{\bf w}\in \overline T~:~h(S_{\bf w})=h_0\} $$
\vskip -.1in\noindent is compact and connected.}\medskip}

\midinsert
\centerline{\qquad\psfig{figure=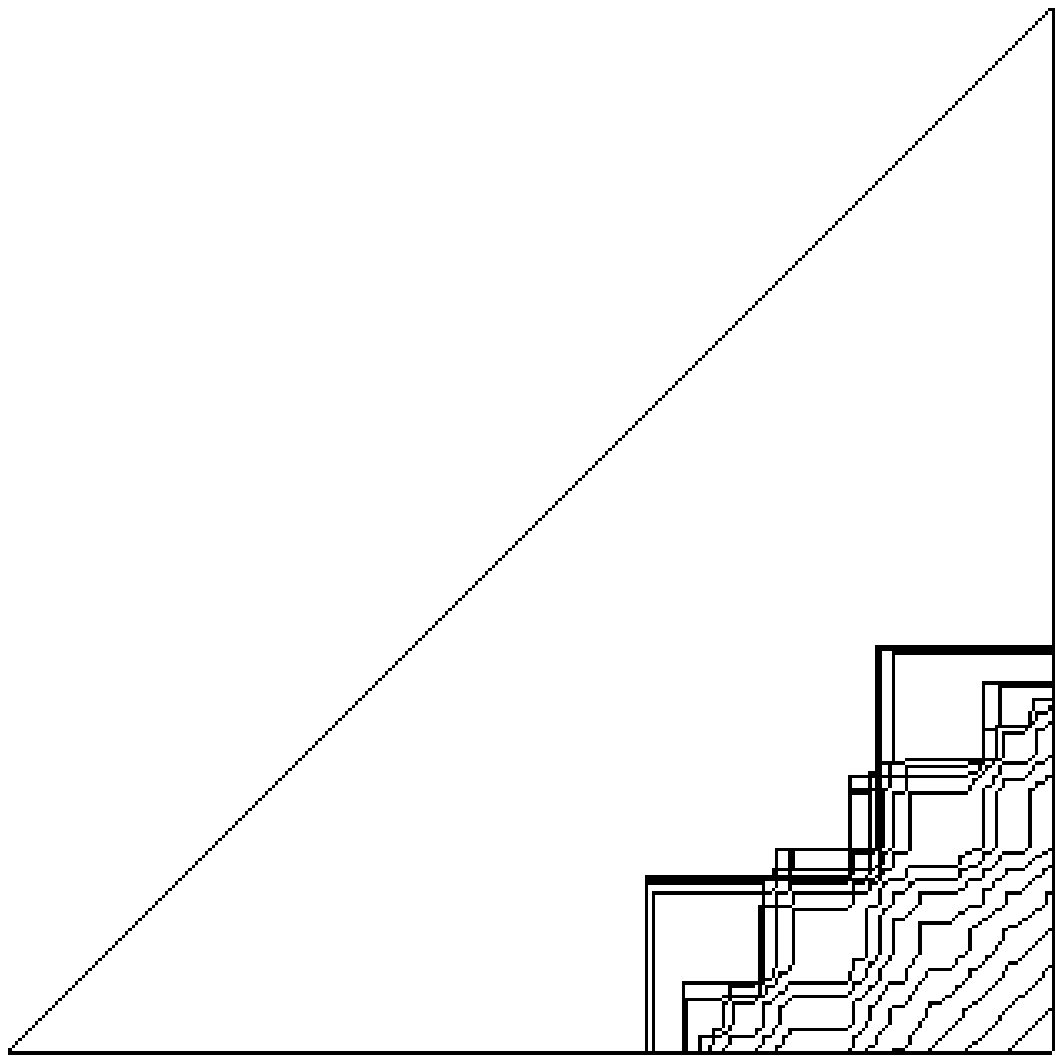,height=2.8in} \hfil
\psfig{figure=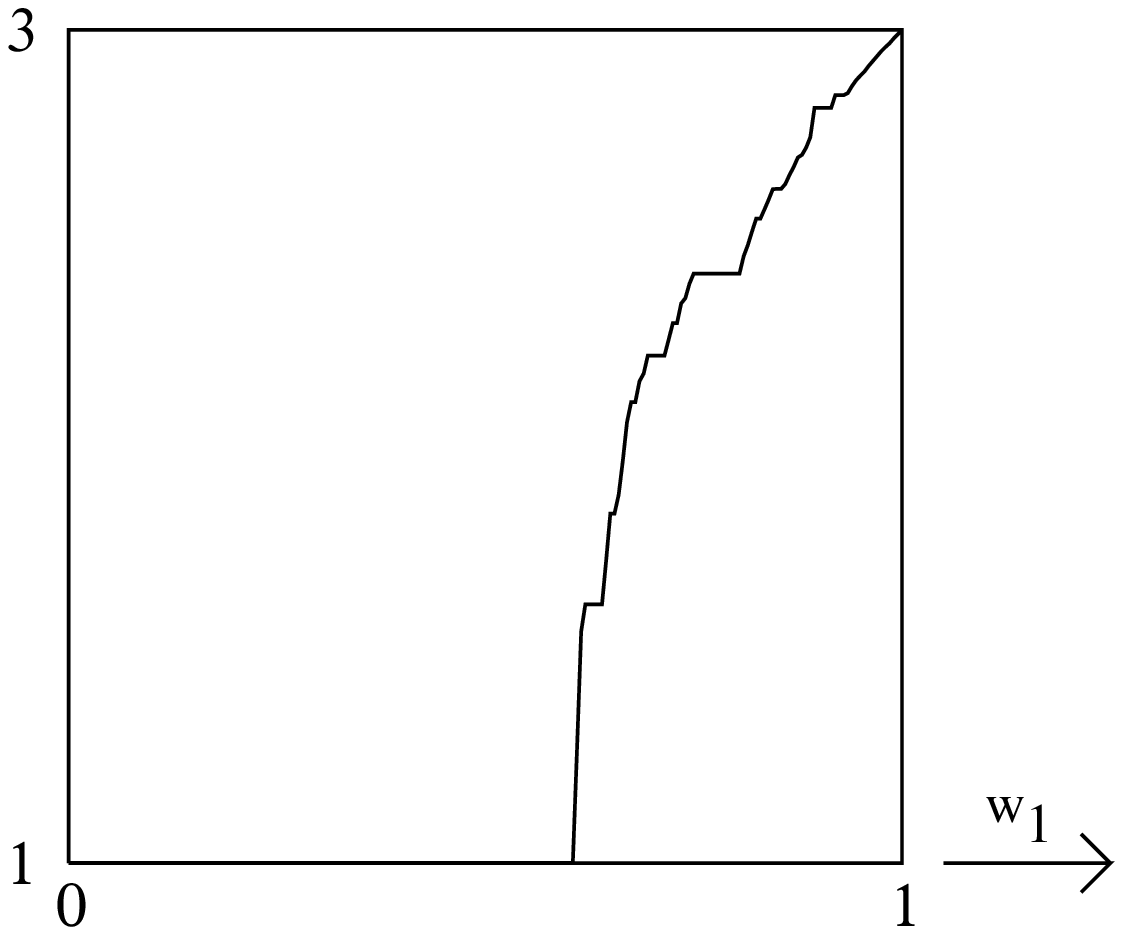,height=2.7in}\qquad}
\centerline{\bit {\bf Left:} Isentropes $\,s={\rm constant}\,$
in the stunted sawtooth parameter triangle, with
$\Delta s=0.1$.}
\centerline{\bit {\bf Right:} Graph of $s=\exp(h)$
as a function of $w_1$ along the bottom edge $w_2=0$ of $\overline T$.}
\endinsert 

{\bf Proof of Theorem 1.} Compactness is clear, since the entropy
function is continuous. Let ${\cal I}^+$ be the union of the
isentrope ${\cal I}(h_0)$ with those seqments of the edges $w_1=1$
or $w_2=0$ on which $h\le h_0$. It follows from (4) that for each line
$w_1+w_2={\rm constant}\;$ which intersects ${\cal I}(h_0)$
the intersection must consist of a point or a closed connected interval.
It follows that we can deformation retract the entire triangle $\overline T$
onto ${\cal I}^+$ by pushing each point towards ${\cal I}^+$ along such
a line $w_1+w_2={\rm constant}$. To check the continuity of this deformation,
it is convenient to rotate the parameter triangle $45^\circ$ by taking
$w_2+w_1$ and $w_2-w_1$ as independent parameters. Then the upper and
lower boundaries of the isentrope will be (not necessarily disjoint)
Lipschitz curves  with $|{\rm slope}|\le 1$, and it follows easily
that our deformation is continuous. Finally, we can certainly
deformation retract ${\cal I}^+$ onto ${\cal I}(h_0)$. Since $\overline T$ is
contractable, it follows that ${\cal I}(h_0)$ is contractable, and hence
connected.\QED\medskip

{\bf Remark.} This statement for the $+-+$
sawtooth family generalizes naturally to $m\!$-modal stunted sawtooth
maps for any $m\ge 1$ and for either one of the two possible shapes.
However, the proof is somewhat harder in the general case. (See [DGMT].)
\eject

\centerline{\bf 5. ``Bones'' in the Parameter Triangle.}\smallskip

The stunted sawtooth family is well understood, but the corresponding
cubic family is poorly understood.
In order to relate theses two families, we
introduce some terminology from MacKay
and Tresser [Ma\s T]. By a {\bit bone} in the parameter triangle $\overline T$
we mean the compact set consisting of all parameter
values for which a specified critical point has periodic orbit with
specified order type. More precisely, the {\bit left bone\/} $B_-(\O)$ is
the set of
parameter values for which the left hand critical point is periodic with
{\bit order type} $\O$.  By definition,
this means that the points of the orbit, numbered as $x_1<\cdots<x_p$,
satisfy $x_i\mapsto x_{\O(i)}$ where $\O$ is some given cyclic permutation of
$\{1,\ldots,p\}$. The {\bit dual right bone\/}
$B_+(\O)$ is the set of parameter values
for which the right critical point is periodic with this same order type.
We will usually assume that the period $p$ is two
or more, and we only allow those order types which can actually occur
for a bimodal map of shape $+-+$.
These definitions make sense either for the stunted sawtooth family or
for the cubic family. (By definition, we take the center points $1/3$
and $2/3$ of the two plateaus as the ``critical points'' for the
stunted sawtooth map.)
We will insert the superscript $\bf saw$ respectively
$\bf cub$ in order to distinguish these two cases.

\midinsert
\centerline{\psfig{figure=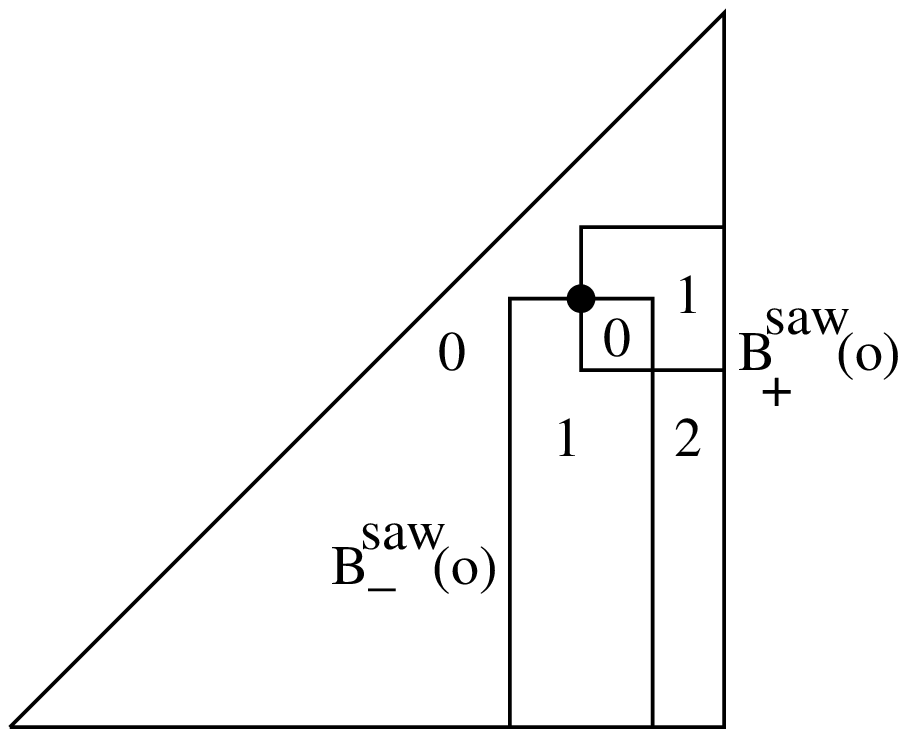,height=2.2in}}
{\QP\bit Dual bones for the stunted sawtooth family. There is a
preferred intersection point, called the common ``center point'' of these
bones, such that the two critical points of the associated map
belong to a common periodic
orbit. It is marked by a heavy dot in the figure.
(The complementary regions have
been labeled by the number of `negative' periodic orbits
with the given order type $\O$. Compare \S8.)\par}
\endinsert

Note that two left bones, or two right bones, are disjoint,
almost by definition.
For the stunted sawtooth family, we have the following simple description:
\smallskip

{\QP{\bf Lemma 3.} {\it Each non-vacuous bone $B_\pm^{\rm saw}(\O)$ of period
$p\ge 2$ is a simple arc with both endpoints on a common edge of the
triangle $\overline T$, and is made up out of three line segments which are
alternately horizontal and vertical. Any pair $B_-^{\rm saw}(\O)$ and
$B_+^{\rm saw}(\O')$ intersect transversally in either 0, 2, or 4 points.
Dual bones always intersect in exactly two points, as illustrated.}\smallskip}

(Compare the schematic picture above.) The proof is not difficult.\QED

{\bf Remark.} Here it is essential that we exclude the case $p=1$,
which behaves quite differently.
(For maps of shape $-+-$ the case $p=2$ is also different, so one must
assume that $p\ge 3$.)  In the $+-+$ case, note that a map
with a critical fixed point must necessarily have zero entropy.

The corresponding statement for the cubic family is much more difficult.
First note the following.

{\QP{\bf Lemma 4.} \it Each bone $B_\pm^{\rm cub}(\O)$ is a smooth
1-dimensional manifold with exactly two boundary points, and these
boundary points belong to a common edge of $\overline T$ (provided that
$p\ge 2$). Any intersection between bones is necessarily transverse.
\smallskip}

{\bf Remark.} Evidently the intersection points of bones are precisely
those points in parameter space for which {\it both\/} critical
points are periodic. In general, the two critical points will belong to
disjoint periodic orbits. The only exception is for the preferred
intersection point of two dual bones. In this exceptional case, the
two critical points belong to a common orbit.

{\bf Proof of Lemma 4.} First
consider the corresponding statement for the family
of complex maps $z\mapsto z^3-3a^2z+b$, with critical points $\pm a$.
It is proved in [M3] that the locus ${\cal S}_\pm(p)$ of points for which
$\pm a$ has period $p$ is a smooth complex curve. Furthermore, for each
$p$ and $q$ the curves ${\cal S}_+(p)$ and ${\cal S}_-(q)$ intersect
transversally. In fact, ${\cal S}_+(p)$ has transverse intersection with
any curve consisting
of points for which the other critical point $-a$ is preperiodic. (The proofs
make essential use of quasi-conformal surgery. Compare [St], where analogous
results for quadratic rational maps are proved by similar methods.)

Restricting to the real $(a,b)\!$-plane, we obtain a corresponding
statement for real cubic maps: {\it In the family of real maps\vskip -.2in
$$x~\mapsto~ x^3-3a^2x+b~,\eqno (6)$$\vskip -.1in
\noindent
the locus of pairs $(a,b)$ for which $a$ (or $-a$) is periodic of period
$p$ forms a smooth $1\!$-dimensional manifold
without boundary.}\bigskip

\centerline{\qquad\qquad\psfig{figure=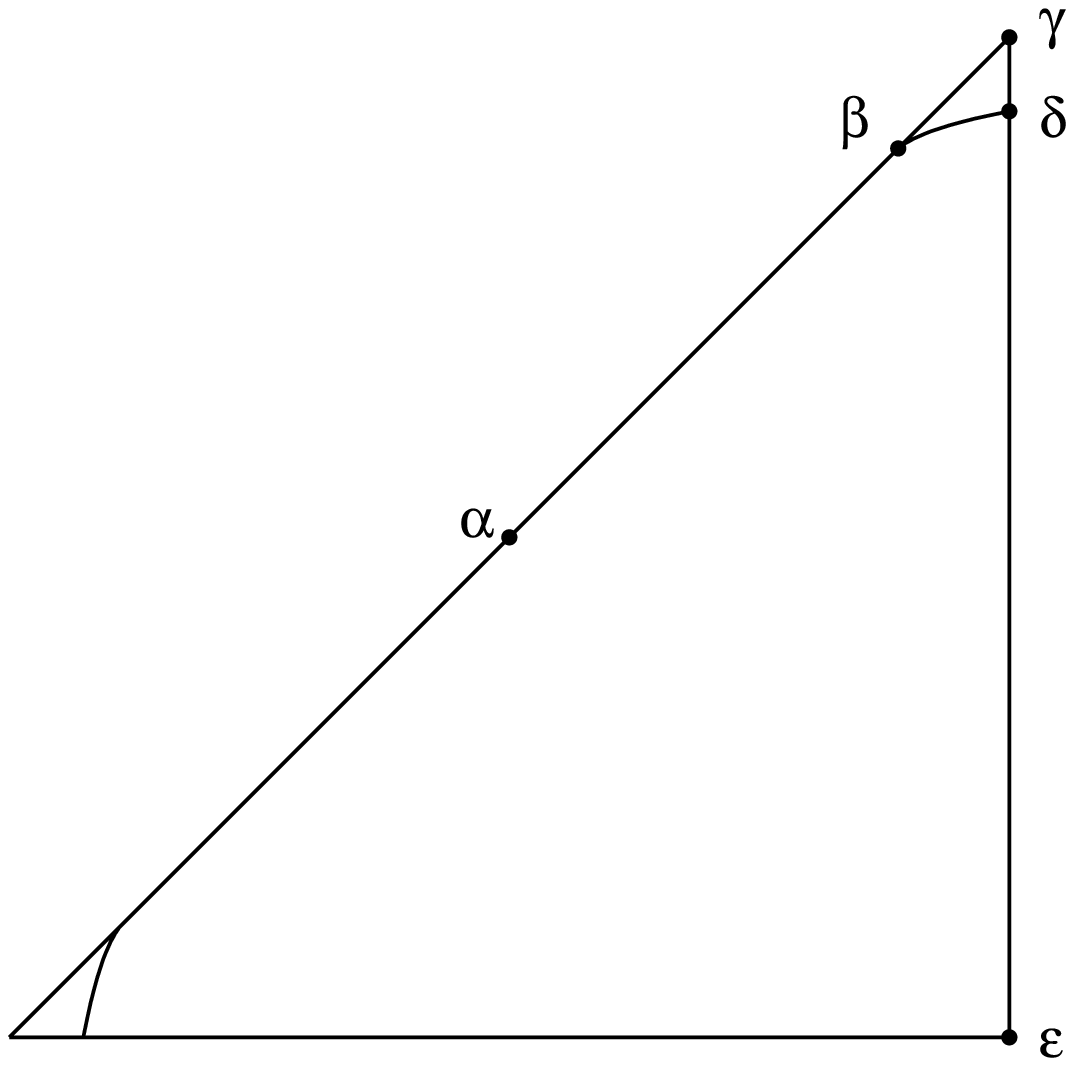,height=2in}
        \hfil
	\psfig{figure=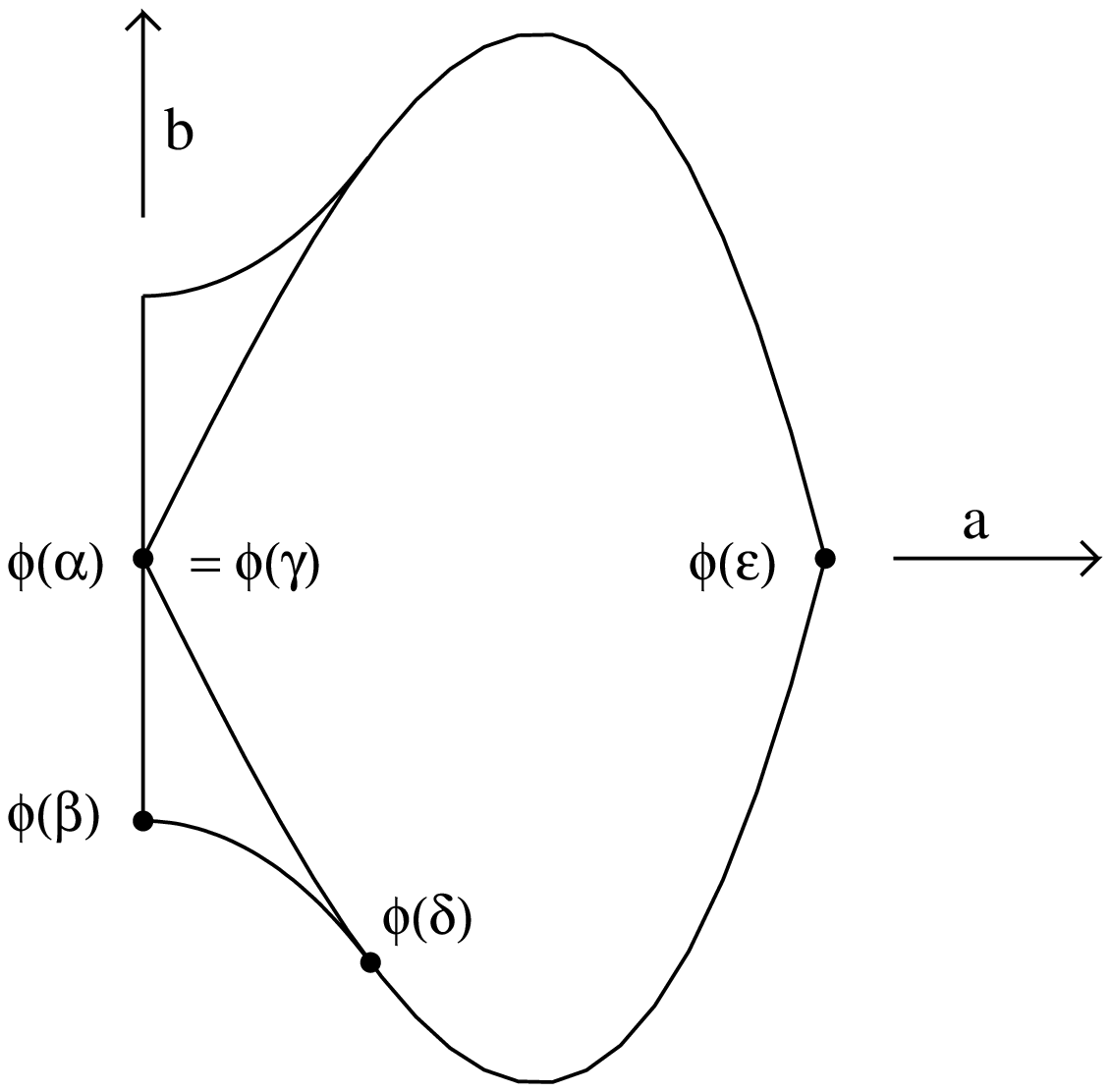,height=2in}\qquad\qquad}

{\QP\bit The parameter triangle ~$\overline T\,$, and its image
~$\phi(T)$~ in the $(a,b)\!$-plane. The two edges $\alpha\beta$ and $\beta
\gamma$ fold together in the negative $b\!$-axis, and the triangle
$\beta\gamma\delta$ folds over so that its image is covered twice.\par}
\eject

In order to relate this to the family of cubic maps $f_\v$,
we note that each $f_\v$ is {\bit positively affinely conjugate\/} to a unique
map in the normal form (6) with $a>0$. That is, there is a unique affine map
$L(x)=c\,x+d$ with $c>0$ so that $L\circ f_\v\circ L^{-1}$ has the
required form. It follows that there is a well defined smooth map
$\phi:T\to\R^2$ which associates to each $\v\in \overline T$ the associated
pair $\phi(\v)=(a,b)$. Evidently
the pre-image of the curve ${\cal S}_\pm(p)$ under $\phi$
is the union of all bones
$B_\pm^{\rm cub}(\O)$ of period $p$. If $\phi$ were a diffeomorphism,
then it would
follow immediately that each bone is a smooth manifold. In fact, the
situation is more complicated, since $\phi$ folds over two of the corners
of the triangle $\overline T$. Thus $\phi$ fails
to be a local diffeomorphism along two fold curves, which correspond
to values $\v$ for which the graph of $f_\v$ is
tangent to the diagonal at one of its two boundary fixed points.
However, for $\v$ along these fold curves or in the folded over regions,
both critical orbits converge to one fixed point, so that
neither critical point can be periodic of period $\ge 2$. It follows easily
that each $B_\pm^{\rm cub}(\O)$
is indeed a smooth 1-manifold, with boundary points at most on the
horizontal or vertical part of the boundary
of $\overline T$.

\midinsert
\centerline{\psfig{figure=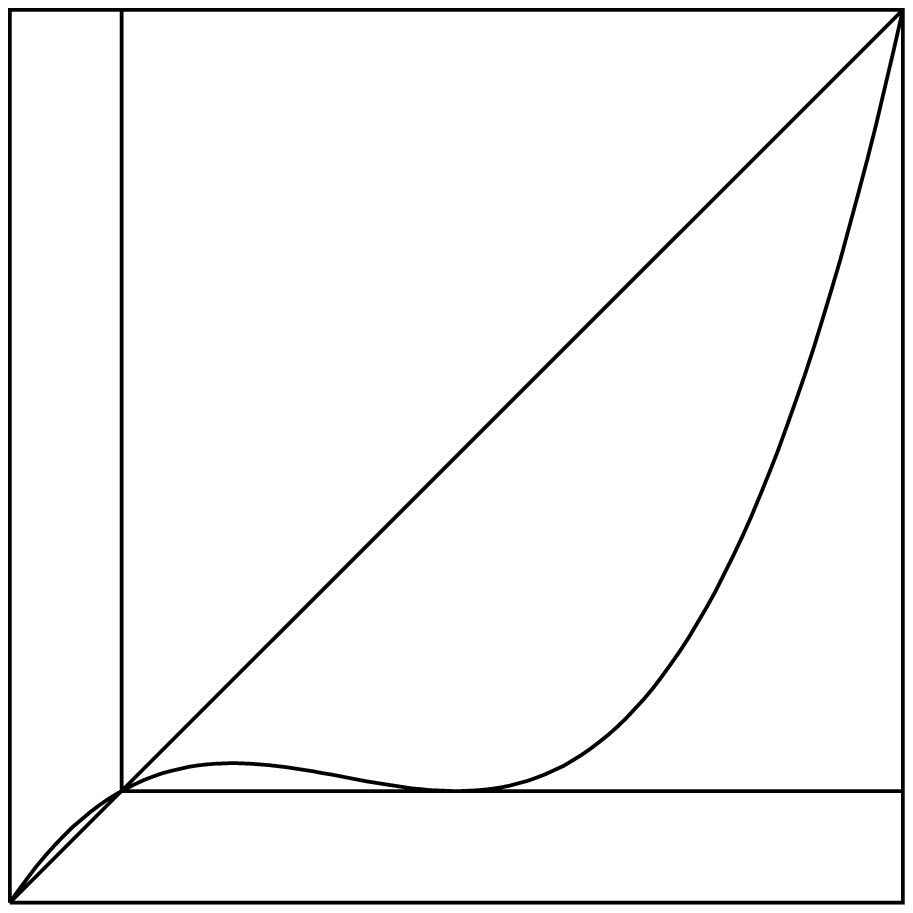,height=2.4in}}
{\QP\bit A cubic map which carries two distinct intervals bimodally
onto themselves. Such a map is positively affinely conjugate
to $f_\v$ for two distinct critical value vectors $\v$. These
two points $\v\in \overline T$ will be folded together by the map $\phi$.
Note that the dynamics of such a map $f_\v$
is necessarily rather trivial.\par}
\endinsert

In order to analyze these possible boundary points $B_\pm^{\rm cub}\cap
\partial \overline T$, we will need to make use of the following basic principle.
By definition, a piecewise monotone map is {\bit post-critically finite\/}
if the orbit of every critical point is periodic or eventually periodic.

{\QP{\bf Thurston's Theorem for Real Polynomial Maps.} {\it Given any
post-critically finite $m\!$-modal map, there exists one and
up to positive affine conjugation only one polynomial of degree $m+1$
with the same kneading data.}\smallskip}

Although this statement is well known to experts, it is difficult
to find in the literature. The proof makes essential use of
complex methods. In fact this statement is an easy corollary of a
much more complicated statement for complex polynomials (or for
complex rational maps). For further information, the reader is referred
to [DH1], [Po] and [dM\s vS], as well as the discussion in [M\s Th].\QED

Assuming Thurston's Theorem, the proof of Lemma 4 concludes as follows.
Consider a map $f_\v$ belonging to the intersection $B_-^{\rm cub}(\O)\cap
\partial \overline T$. Thus the left critical point of $f_\v$ is periodic,
with some period $p\ge 2$. It is easy to check that the two critical points
cannot coincide. (Here we make use of the fact that our maps have shape $+-+$.
For maps of shape $-+-$ we would rather have to assume that $p\ge 3$ in order
to avoid the case $c_1=c_2$.) Since $f_\v\in\partial \overline T$, it follows
that either $v_1=1$ or $v_2=0$. In other words, at least
one of the two critical points must map to a boundary fixed point.
But $c_1$ is periodic, so it follows that we must be on the edge $v_2=0$.
If we cut the interval $I$ at the points of the period $p$ orbit $\{c_1,v_1,
f_\v(v_1),\ldots\}$, then since our map must have shape $+-+$ a
little work shows that there are exactly two of the $p+1$
complementary intervals
where the right hand critical point can be placed. Using Thurston's Theorem,
it follows that the intersection $B^{\rm cub}_-(\O)\cap\partial \overline T$
consists of exactly two points. Further, this intersection is
transverse, so these two intersection points must belong to the
boundary $\partial B^{\rm cub}_-(\O)$.
\QED\medskip

\centerline{\bf 6. Three Conjectures.}\smallskip

It follows from Lemma 4 that each bone $B^{\rm cub}_\pm(\O)$
consists of a simple arc, possibly
together with one or more disjoint simple closed curves, which we
may call ``bone-loops''. Thus a hypothetical {\bit bone-loop\/}
would be a simple closed curve in the parameter triangle $\overline T$
consisting entirely of points $\v$ for which one critical point of $f_\v$ (say
the left one) is periodic. In fact we conjecture
that such bone-loops do not exist:

{\QP
{\bf Connected Bone Conjecture.} \it Every bone $B_\pm^{\rm cub}(\O)$ for the
cubic family is a simple arc.\smallskip}

Although we cannot prove this conjecture, we will show that it would
follow from a well known classical conjecture.\smallskip

{\bf Generic Hyperbolicity.} Recall that a polynomial or
rational map is said to be {\bit hyperbolic\/} if the orbit of every complex
critical point converges towards an attracting periodic orbit. By the
{\bit Generic Hyperbolicity Conjecture\/} for some given family $\cal F$
of maps we mean
the conjecture that every map $f_0\in{\cal F}$ can be approximated arbitrarily
closely by a map $f\in{\cal F}$ which is hyperbolic. (Compare the discussion
in [F, p. 73].) For the family of real quadratic maps,
a proof of this statement has been announced by Swiatek [S], and more recently
by Lyubich [L]. (See also [Mc2], which proves a weaker but closely related
result.) However, we will rather need the statement for cubic maps. In fact,
the main result of this note will be to relate
these conjectures to the Monotonicity Conjecture, as stated in \S3:\medskip

\centerline{\bit Generic Hyperbolicity for Real Cubic Maps}
\smallskip

\centerline{$\Longrightarrow$ \bit Connected Bone Conjecture}\smallskip

\centerline{$\qquad\Longrightarrow$ \bit Monotonicity Conjecture.}
\medskip

\noindent
The proofs of these implications will be given in \S7 and \S8 respectively.
In order to carry out these proofs, we must first develop a closer relationship
between the cubic family and the stunted sawtooth family.
\bigskip\medskip

\centerline{\bf 7. Intersections of Bones, and the $n\!$-Skeleton.}\smallskip

By the $n\!$-{\bit skeleton\/} $S^{\rm saw}_n$ for the stunted
sawtooth family, we will mean the union of all bones $B^{\rm saw}_\pm(\O)$
of period $2\le p\le n$, together with the boundary $\partial \overline T$. In the
case of the cubic family it will be convenient to modify this definition
slightly as follows. Recall that each bone $B^{\rm cub}_\pm(\O)$ consists
of a simple arc, possibly together with some disjoint
simple closed curves. We will use
the notation $A_\pm(\O)$ for this simple arc. By the $n\!$-{\bit skeleton\/}
$S^{\rm cub}_n$ we will mean the union of all of these simple arcs $A_\pm(\O)$
with period $2\le p\le n$, together with $\partial \overline T$. For an
analysis of the structure of these skeletons, see [R\s S], [R\s T],
as well as [Ma\s T]. We will prove the following.

{\QP{\bf Theorem 2.} \it For each $n\ge 2$ there exists
a homeomorphism $\eta_n$ from the parameter
triangle for stunted sawtooth maps onto the parameter triangle for cubic
maps which carries each bone $B_\pm^{\rm saw}(\O)$ of period $2\le p\le n$
onto the corresponding bone-arc $A_\pm^{\rm cub}(\O)$, and hence carries
the skeleton $S_n^{\rm saw}$ homeomorphically onto $S_n^{\rm cub}$.
\par}
\midinsert
\centerline{\psfig{figure=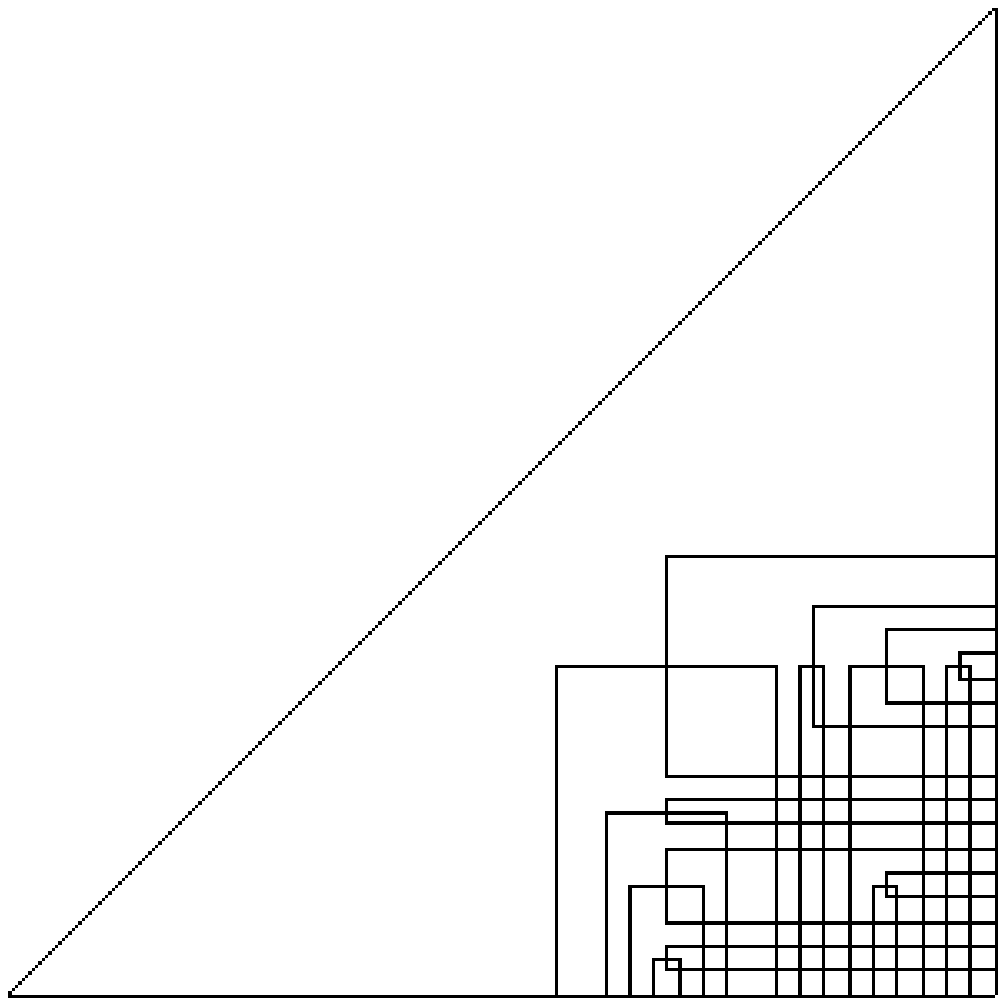,height=2.7in}
\psfig{figure=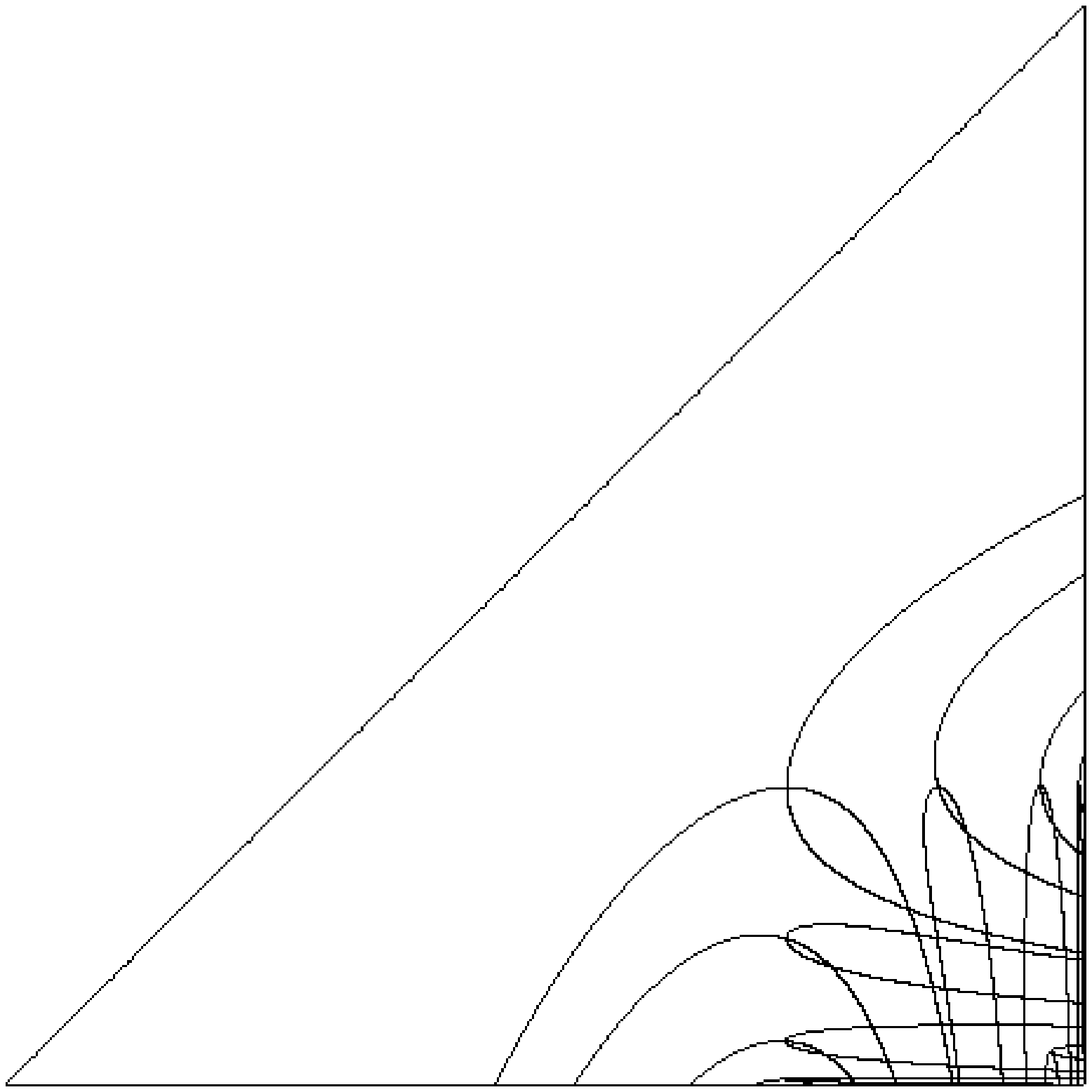,height=2.7in}}
{\QP\bit Bones of period $2, 3, 4$ for the stunted sawtooth and cubic
families. (As we traverse the bottom edge of $\overline T$ from left to right,
we meet bones of period $2,4,3,4,4,3,4,2,4,4,3,4,4,3,4,4$ respectively.)
In the stunted sawtooth family, the dual bone can be recognized as the
unique bone which meets the middle edge of a given bone.\par}
\vskip -.1in\endinsert

{\bf Outline
Proof.} By a {\bit vertex\/} of the skeleton $S_n^{\rm cub}$
or $S_n^{\rm saw}$ we will mean either an endpoint of a bone or a point
of intersection between two
bones (necessarily one left bone and one right bone). Evidently each such
vertex corresponds to a post-critically finite map. Every such map is
uniquely determined by its kneading data. In fact, any kneading data which
can occur for an arbitrary post-critically finite $+-+$
bimodal map must actually occur for one and only one map in this family.
In the cubic case, this follows easily
from Thurston's Theorem as stated in \S5, while in
the stunted sawtooth case it follows from an elementary argument
as in the proof of Lemma 2. Thus there is a natural one-to-one
mapping from the vertices of the skeleton $S_n^{\rm saw}$ onto
the vertices of $S_n^{\rm cub}$.

Remark: Here it is essential that we only consider bones of period $\ge 2$,
so that none of these vertices lie on the upper left edge $v_1=v_2$
of the parameter triangle where things
behave somewhat strangely. (Compare the proof of Lemma 6 below.)

Next we must check that these vertices
occur in the same order as we traverse some bone $B_\pm^{\rm saw}(\O)$
or as we traverse the corresponding bone-arc $A_\pm(\O)$ from one
endpoint to the other. By the {\bit center point\/}
of a bone we will mean the unique
point for which both critical points belong to a common orbit. (This is
one of the two intersection points with the dual bone, and is evidently
also the center point for the dual bone.) We will see that the center
point of any bone $B^{\rm cub}_\pm(\O)$ necessarily belongs to the
component $A_\pm(\O)$. Now, as we traverse any $B^{\rm saw}_\pm(\O)$
or $A_\pm(\O)$ from one end to the other, we claim that the entropy
decreases monotonically until we reach the center point, and then
increases monotonically until we reach the other end. In other words,
if we divide each $B^{\rm saw}_\pm(\O)$
or $A_\pm(\O)$ into two {\bit half-bones\/}
by cutting at this center point, then along each half-bone
the entropy changes monotonically.
In fact, the kneading
sequence for the non-periodic critical point changes monotonically
along each half-bone. For the stunted sawtooth family,
this can be proved by a direct argument. For the cubic family, it
can be proved by an argument which is completely analogous to the proof
in [M\s Th] of monotonicity for the quadratic family. Note first
that the kneading sequence for the periodic critical point is fixed
as we traverse the half-bone (although it is different from one half-bone
to the other). The kneading sequence for the remaining
critical point varies continuously with the cubic
map, except for discontinuities of a very special form at those maps for which
this remaining critical point eventually maps to one or the other critical
point. Using these facts, we obtain an ``intermediate value theorem''
for admissible kneading sequences as we traverse the half-bone. In particular,
any post-critically finite kneading data which can occur for any $+-+$
bimodal map must occur somewhere along the appropriate half-bone.
It now follows
that the kneading sequence must change monotonically. For otherwise
some post-critically finite kneading sequence would have to occur twice,
which is impossible by Thurston's Theorem.
The rest of the proof is reasonably straightforward.\QED

As an immediate corollary to this argument we obtain the following.

{\QP{\bf Lemma 5.} \it A bone-loop in the cubic parameter triangle
cannot contain any post-critically finite point. In particular, it
cannot intersect any other bone.\smallskip}

{\bf Proof.}  We have shown in the proof of Theorem 2
that all possible kneading types of post-critically finite points in a bone
$B_\pm^{\rm cub}(\O)$ can be found somewhere
along the bone-arc $A_\pm(\O)$. There cannot be any other post-critically
finite points by Thurston's Theorem.\QED

{\QP{\bf Lemma 6.} \it The region enclosed by a bone-loop in $\overline T$
cannot contain any hyperbolic maps.\smallskip}

{\bf Proof.} By a {\bit hyperbolic component\/} $H$ in the cubic
parameter triangle $\overline T$
we mean a connected component in the open set consisting of all
$\v\in \overline T$ such that $f_\v$ is hyperbolic. It follows immediately
from [M2] that every hyperbolic component which lies completely
in the interior of $\overline T$ is an
open topological 2-cell which contains a unique post-critically
finite point called its {\bit center\/}. (Compare [Mc1].)
Furthermore, every bone which
intersects such a hyperbolic component $H$ must intersect it in a
simple arc passing through the center point of $H$. (Compare [M3].
There may be just one bone passing through the component, or there may
be two which intersect transversally at the center point.)

{\bf Note:} Here we exclude the three exceptional hyperbolic components which
meet the boundary of $\overline T$. These are ``centered'' at the mid-point
and end-points of the
edge $v_1=v_2$. Because of the folding mentioned in the proof of
Lemma 4, these three distinct center points actually correspond to
just one affine conjugacy class of maps, which is represented
for example by $x\mapsto x^3$.

Suppose now that some hyperbolic component $H$ intersects
the region enclosed by a bone-loop $L$. Then $H$ certainly cannot be one
of the three exceptional hyperbolic components which meet the boundary
of $\overline T$. If $H$ intersects the loop
itself, then its center point must belong to $L$, which is impossible by
Lemma 5. On the other hand, if $H$ is completely enclosed by $L$ then
we can choose a bone-arc $A_\pm(\O)$ which passes through the
center point of $H$. By the Jordan Curve Theorem, $A_\pm(\O)$ must
intersect L, which again contradicts Lemma 5.\QED

{\QP\bf Theorem 3. {\it The Generic Hyperbolicity Conjecture
for real cubic maps\break implies the Connected Bone Conjecture.}\smallskip}

{\bf Proof.} This statement, which was promised in \S6, clearly
follows as an immediate corollary to Lemma 6.\QED

{\bf Remark.} There is real hope that something much less than the full
Generic Hyperbolicity Conjecture for cubic maps
might be enough to prove the Connected Bone Conjecture.
In particular, it would be enough to prove generic hyperbolicity along
each bone. Thus we need only study very special one-parameter families of
cubic maps.
The known techniques for dealing with quadratic maps might well suffice
to deal with this special case.
\bigskip

\centerline{\bf 8. Bones, Negative Orbits, and the Monotonicity Conjecture.}
\smallskip

In order to relate bones and entropy, let us first recall the following
result from [M\s Th]. A fixed point $f(x_0)=x_0$ for
a piecewise monotone
map will be called either {\bit positive\/}, {\bit negative\/}, or {\bit
critical\/} according as the map $f$ is monotone increasing or monotone
decreasing throughout a neighborhood of $x_0\,$, or has a turning point at
$x_0$. {\bf Definition:} Let $N(f)$ be the number
of critical fixed points of $f$ plus
twice the number of negative fixed points. This number is closely related
to the total number of fixed points $\#{\rm fix}(f)$, but
is easier to work with since
it is more robust. Note that $N(f)$ is always finite, in fact
$N(f)\le \ell(f)+1$ where $\ell(f)$ is the number
of laps. It is easy to check that
the sequence of numbers $N(f^{\circ k})$ is completely
determined by the kneading data for $f$, as described in \S4. We will need the
following estimate, which is similar to formula (3) of \S3
but true in much greater generality.

{\QP{\bf Lemma 7.} \it The topological entropy of an arbitrary
piecewise monotone map is given by the formula\vskip -.2in
$$	h(f)~=~\limsup_{k\to\infty} {\log N(f^{\circ k})\over k}~. $$
\smallskip}

{\bf Proof.} This follows easily from [M\s Th], or from [Pr] or [B\s R].\QED

{\QP{\bf Lemma 8.} \it 
Let $\v$ and $\v'$ be two points in the cubic parameter triangle
such that the associated maps have entropy $h(f_\v)\ne h(f_{\v'})$.
Then any path from $\v$ to $\v'$ in the parameter triangle
$\overline T$ must cross infinitely many bones.\smallskip}

For according to Lemma 7 the difference $|N(f_\v^{\circ k})-
N(f_{\v'}^{\circ k})|$ must be unbounded as $k\to\infty$. But clearly,
as we deform $\v$ along some path in $\overline T$, the number
$N(f_\v^{\circ k})$, which measures the number of decreasing laps of
$f_\v^{\circ k}$ whose graph crosses the diagonal, will remain constant
except as we pass through a map which has a critical periodic orbit
of period dividing $k$. In other words, $N(f^{\circ k})$ remains
constant unless we pass through a bone of period dividing $k$. (Here
bones of period 1 must also be allowed, but cause no difficulty.)
Further details are easily supplied.\QED

Completely analogous arguments apply to the stunted sawtooth family,
provided that we define the concept of a ``negative'' fixed point
of $S_\w^{\circ k}$ in an appropriate formal manner. {\bf Definition:}
A fixed point $S_\w: x_0\mapsto x_1\mapsto\cdots\mapsto x_k=x_0$ of
$S_\w^{\circ k}$
is {\bit critical\/} if the orbit passes through one of the critical
points $1/3$ or $2/3$, and otherwise is {\bit positive\/} or {\bit negative\/}
according as the number of $x_i$ in the interval $(1/3\,,\,2/3)$
is even or odd. Using this definition, the analogues of Lemmas 7 and 8 are
easily verified.

Associated with the skeleton $S_n^{\rm saw}\subset \overline T$ is
a {\bit topological cell structure\/} on $\overline T$. That is, we can
partition $\overline T$ into subsets, each of which is homeomorphic
either to a
point, an open interval, or a 2-dimensional open unit disk. Furthermore,
these subsets
fit together nicely so that the closure of each one is topologically
a point, a closed interval, or a closed 2-disk.
By definition, the open 2-cells in this cell structure are the
connected components of the complement $\overline T\smallsetminus
S_n^{\rm saw}$, the 0-cells are the vertices as described in \S7,
and the open 1-cells are the connected components of $S_n^{\rm saw}
\smallsetminus\{{\rm vertices}\}$. The resulting cell complex will be denoted
by $\overline T_n^{\,\rm saw}$. There is a completely analogous cell
complex $\overline T_n^{\,\rm cub}$ for the cubic family. These two cell
complexes are homeomorphic by Theorem 2, in a homeomorphism which
takes each vertex to a vertex of the same entropy and each edge
to an edge with the same interval of entropies.

{\QP{\bf Lemma 9.} \it For each $\epsilon>0$ there exists an integer
$n$ so that, for any two points $\w$ and $\w'$ belonging to a common
closed cell of the complex $\overline T_n^{\,\rm saw}$, we have
$$	|h(S_{\w'})-h(S_\w)|~<~\epsilon~. $$
If the Connected Bone Conjecture is true, then there is an
analogous statement for the complex $\overline T_n^{\,\rm cub}$ and the
family of cubic maps $f_\v$.\smallskip}

{\bf Proof.} Otherwise we could find $\epsilon>0$ so that for each $n$
there existed points $\w_n$ and $\w'_n$ in a common cell of $\overline
T_n^{\,\rm saw}$ with $|h(S_{\w'_n})-h(S_{\w_n})|\ge \epsilon$. After passing
to infinite subsequences, we could assume that both sequences converge,
say $\w_k\to\w$ and $\w'_k\to\w'$,
and furthermore that all of the $\w_k$ and $\w'_k$
belong to a common cell of $\overline T_n^{\,\rm saw}$ whenever $k\ge n$.
Thus the limit points $\w$ and $\w'$ would belong to a common closed
cell of $\overline T_n^{\,\rm saw}$ for every $n$, but by continuity
the associated entropies would differ by at least $\epsilon$. This is
impossible by Lemma 8. The proof for the cubic family is similar.\QED

{\QP{\bf Lemma 10.} \it The entropy function $\w\mapsto h(S_\w)$ for the
stunted sawtooth family, restricted to any closed cell of the cell
complex $\overline T_n^{\,\rm saw}\!$, takes its maximum
and minimum values on the boundary (and in fact on the set
of vertices). If the Connected Bone Conjecture is true, then
the analogous statement holds for the entropy function $\v\mapsto h(f_\v)$
for cubic maps, restricted to any
closed cell of the cell complex $\overline T_n^{\,\rm cub}\!$.
\smallskip}

{\bf Proof.} For the stunted sawtooth family, this
follows easily from Theorem 1, together with
the fact that the entropy function is monotone along each edge. (Compare
the proof of Theorem 2.) To prove the analogous statement for the cubic
family, suppose for example
that for some point $\v_0$ of a closed cell $C$, the value
$h(f_{\v_0})$ were strictly larger than the maximum value $h_{\rm max}$
on the boundary of $C$. Let $2\epsilon = h(f_{\v_0})-h_{\rm max}$.
According to Lemma 9 we can choose $m>n$ so that $h$ varies by less than
$\epsilon$ on each cell of $\overline T_m^{\rm cub}$. Let $C'\subset C$
be a cell of $\overline T_m^{\rm cub}$ which contains this point $\v_0$,
and let  $\v'$ be any vertex of $C'$. Then it follows that
$h(f_{\v'}) > h_{\rm max}$.
Since the homeomorphism $\eta_m$ of Theorem 2 carries
vertices to vertices with the same entropy, this would yield a vertex in
the complex $\overline T_m^{\rm saw}$ satisfying a corresponding
inequality. But this is impossible,
since Lemma 10 is known to be true for the stunted sawtooth family.\QED

{\QP{\bf Theorem 4.} \it If the Connected Bone Conjecture is true, then
every isentrope $\{\v\in T : h(f_\v)=h_0\}$ for the cubic family
is connected. Thus, the Connected Bone Conjecture implies the
Monotonicity Conjecture.\smallskip}

{\bf Proof.} For each $n$, the union of all closed cells of $\overline
T_n^{\,\rm saw}$
which touch the $h_0\!$-isentrope forms a compact set, which is connected
by Theorem 1. The
corresponding union of cells in $\overline T_n^{\,\rm cub}$ forms a compact
connected set by Theorem 2, and
contains the corresponding isentrope by Lemma 10.
The intersection of these sets, as $n\to\infty$, will be precisely equal
to the required isentrope by Lemma 9. Since an intersection of compact
connected sets is compact and connected, the conclusion follows.\QED

\medskip

\centerline{\bf Appendix on Computation.}\smallskip

As a supplement to this paper, the computer programs which were used to
make figures are available, in documented form, via ftp. They include
subroutines to compute the topological entropy of an arbitrary unimodal
or bimodal map. (The latter is based on Block and Keesling [B\s K].)
For further information, send email to
`IMS@math.sunysb.edu'.\bigskip


\centerline{\bf References.}\smallskip

\ref [B] R. Bowen, {\sl On Axiom A Diffeomorphisms\/}, Proc. Reg. Conf.
Math. {\bf 35}, 1978.

\ref [B\s K] L. Block and J. Keesling, ``Computing
topological entropy of maps of the
interval with three monotone pieces'', {\sl J. Statist. Phys. \bf 66} (1992)
755-774.

\ref [BMT] K. Brucks, M. Misiurewicz, and C. Tresser,
``Monotonicity properties of the family of trapezoidal maps,'' {\sl
Commun. Math. Phys. \bf 137} (1991) 1--12.

\ref [B\s R] V. Baladi and D. Ruelle, ``An extension of the
theorem of Milnor and Thurston on the zeta functions of interval maps'',
{\sl Ergodic Theory and Dynamical Systems}, to appear.

\ref [BST] N.J.Balmforth, E.A.Spiegel and C. Tresser, ``The topological
entropy of one-dimen\-sional maps: approximation and bounds'', to appear.

\ref [D] A. Douady, ``Entropy of unimodal maps'', Proceedings NATO Institute
on Real and Complex Dynamical Systems, Hiller{\o}d 1993, these proceedings.

\ref [D\s G] S.P. Dawson and C. Grebogi, ``Cubic
maps as models of two-dimensional antimonotonicity,'' {\sl Chaos,
Solitons \& Fractals \bf 1} (1991), 137--144.

\ref [DGK] S.P. Dawson, C. Grebogi and  H. Ko{\c c}ak,
``A geometric mechanism for antimonotonicity in scalar maps with two
critical points,'' {\sl Phys. Rev. E} (in press).

\ref [DGKKY] S.P. Dawson, C. Grebogi, I. Kan, H. Ko\c cak and  J.A. Yorke,
``Antimonotonicity: inevitable reversals of period doubling
cascades,'' {\sl Phys. Lett.A} {\bf 162} (1992) 249--254.

\ref [DGMT] S.P. Dawson, R. Galeeva, J. Milnor and C.
Tresser, ``Monotonicity and antimonotonicity for bimodal maps'',
manuscript in preparation.

\ref [DH1] A.Douady and J.Hubbard, ``A proof of Thurston's Topological
Characterization of Rational Maps''; Preprint, Institute Mittag-Leffler 1984. 

\ref [DH2] A. Douady and J. H. Hubbard,``Etude dynamique des
polyn\^omes quadratiques complexes,'' {\sl I} (1984) \& {\sl II}
(1985), {\sl Publ. Mat. d'Orsay}.

\ref [dM\s vS] W. De Melo and S. Van Strien, {\sl One Dimensional Dynamics},
Springer V., 1993.

\ref [F] P. Fatou, ``Sur les \'equations fonctionnelles, II'', {\sl
Bull. Soc. Math. France \bf 48} (1920) 33-94.

\ref [Ga] R. Galeeva, ``Kneading sequences for piecewise linear bimodal
maps'', to appear.

\ref [Gu] J. Guckenheimer, ``Dynamical Systems'', C.I.M.E. Lectures
(J. Guckenheimer,
J. Moser and S. Newhouse), Birkhauser, Progress in Mathematics {\bf 8}, 1980.

\ref [K] A. Katok, ``Lyapunov exponents, entropy and periodic orbits of
diffeomorphisms'', {\sl Pub. Math. IHES \bf 51} (1980) 137-173.

\ref [L] M. Lyubich, ``Geometry of quadratic polynomials: moduli, rigidity,
and local connectivity'', Stony Brook I.M.S. Preprint 1993/9.

\ref [M1] J. Milnor, ``Remarks on iterated cubic maps'', {\sl
Experimental Math. \bf 1} (1992) 5--24.

\ref [M2] ---------, ``Hyperbolic components in Spaces of Polynomial
Maps (with an appendix by A.~Poirier)'', Stony Brook I.M.S.
Preprint 1992\#3.

\ref [M3] ---------, ``On
cubic polynomials with periodic critical point'', in preparation.

\ref [Ma\s T]
R.S MacKay and C. Tresser,``Boundary of topological chaos for bimodal maps of
the  interval,'' {\sl J. London Math. Soc. \bf 37} (1988), 164--81;
``Some flesh on the skeleton: the bifurcation
structure of bimodal maps'', {\sl Physica \bf 27D} (1987) 412-422.

\ref [Mc1] C. McMullen, ``Automorphisms of rational maps'', pp. 31-60 of
{\sl Holomorphic Functions and Moduli I\/}, ed. Drasin, Earle, Gehring, Kra \&
Marden; MSRI Publ.$\,10$, Springer 1988.

\ref [Mc2] C. McMullen, ``Complex dynamics and renormalization'',
preprint, U.C. Berkeley 1993.

\ref [Mis] M. Misiurewicz, ``On non-continuity of topological entropy'',
{\sl Bull. Ac Pol. Sci., Ser. Sci. Math. Astr. Phys.\bf 19} (1971)
319-320.

\ref [M\s Sz] M. Misiurewicz and W. Szlenk, ``Entropy of piecewise monotone
mappings'', {\sl Studia Math. \bf 67} (1980) 45--63 (Short version:
{\sl Ast{\'e}risque \bf 50} (1977) 299--310).

\ref [M\s Th] J. Milnor and W. Thurston, ``On iterated maps of the
interval,'' {\sl Springer Lecture Notes \bf 1342} (1988), 465--563.

\ref [N] S. Newhouse, ``Continuity properties of entropy'', {\sl Ann. Math.
\bf 129} (1989) 215-235 and {\bf 131} (1990) 409-410.

\ref [Po] A. Poirier, ``On post critically finite polynomials, Part II,
Hubbard Trees'', Stony Brook I.M.S. Preprint 1993/7.

\ref [Pr] C. Preston, ``What you need to know to knead'', {\sl Advances Math.
\bf 78} (1989) 192-252.

\ref [Ro] J. Rothschild, ``On the computation of topological entropy'',
Thesis, CUNY 1971.

\ref [R\s S] J. Ringland and M. Schell, ``Genealogy and bifurcation
skeleton for cycles of the iterated two-extremum map of the interval'',
{\sl SIAM J. Math. Anal. \bf 22} (1991) 1354-1371.

\ref [R\s T] J. Ringland and C. Tresser, ``A genealogy for finite kneading
sequences of bimodal maps of the interval'', preprint, IBM 1993.

\ref [St] J. Stimson, ``Degree two rational maps with a periodic critical
point'', Thesis, Univ. Liverpool 1993.


\ref [Sw] G. Swiatek, ``Hyperbolicity is dense in the real quadratic family'',
Stony Brook I.M.S. Preprint 1992/10.

\ref [Y] Y. Yomdin, ``Volume growth and entropy'', {\sl Isr. J. Math. \bf
57} (1987) 285-300. (See also ``$C^k\!$-resolution of semialgebraic
mappings'', ibid. pp. 301-317.)
\vfil

\def\in{\hskip 2.8in}
\def\out{\hfil\break\vskip -.25in}

\in Silvina P. Dawson\out
\in C.N.L.S.\out
\in Los Alamos National Laboratory\out
\in Los Alamos NM 87545\out\bigskip

\in Roza Galeeva\out
\in Mathematics Dept.\out
\in Northwestern Univ.\out
\in Evanston IL 60208-2730\out\bigskip

\in John Milnor\out
\in Institute for Mathematical Sciences\out
\in SUNY\out
\in Stony Brook NY 11794-3651\out\bigskip

\in Charles Tresser\out
\in I.B.M.\out
\in P.O. Box 218\out
\in Yorktown Heights NY 10598\out
\vfil
\eject
\end